\documentclass[10pt,twoside]{article}
\usepackage{amsmath,amsthm}
\usepackage[latin1]{inputenc}
\usepackage{picinpar}
\usepackage{longtable}
\usepackage{amsrefs}
\usepackage[dvips]{graphicx}
\usepackage[mathscr]{eucal}
\usepackage{calc}
\usepackage{ifthen}
\usepackage{enumerate}
\usepackage{amssymb,latexsym}
\usepackage{amsthm}
\usepackage[dvips]{graphics}
\usepackage{color}
\usepackage{tabularx}
\usepackage{makeidx}
\theoremstyle{definition}

 \newtheorem{teor}{Theorem}

\begin{document}
\bigskip

\thispagestyle{plain}
\par\bigskip
\begin{centering}

\textbf{On a generalization of $P_{3}(n)$}

\end{centering}\par\bigskip
 \begin{centering}
 \footnotesize{AGUSTIN MORENO CA\~NADAS}\par\bigskip
 \end{centering}

$\centerline{\textit{\small{Departamento de Matem\'aticas, Universidad Nacional de Colombia}}}$
\centerline{\textit{\small{Bogot\'a-Colombia}}} $\centerline{\textit{\small{amorenoca@unal.edu.co}}}$

\bigskip
\bigskip \small{We present a formula for $P^{s}_{3}(n)$ the number of partitions of a positive number $n$ into $3$ $s$-gonal numbers,
by using representations of posets over $\mathbb{N}$.}


\bigskip

\textbf{1. Introduction}\par\bigskip

In [10] Lehmer denoted $P_{k}(n)$ the number of partitions of a natural number $n$ into $k$ integral squares
$\geq0$, and solved almost completely the equation $P_{k}(n)=1$. Lehmer claimed that the general problem of
finding a formula for $P_{k}(n)$ was a problem of great complexity. The case $k=3$ was studied by Grosswald, A.
Calloway, and J. Calloway in [4], and Grosswald solved (essentially) the problem, giving the number of partitions
of an arbitrary integer $n$ into $k$ squares (taking into account that, he didn't distinguish between partitions
that contains zeros and those that do not) [5].\par\bigskip

The main goal of this paper is to give a formula for $P^{s}_{3}(n)$, where the $s$-th polygonal number of order or
rank $r$, $p^{s}_{r}$ is given by the formula (often, $0$ is included as a polygonal number [3,7]),
\par\bigskip

\begin{centering}
$p^{s}_{r}=\frac{1}{2}[(s-2)r^{2}-(s-4)r]$.\par\bigskip

\end{centering}

We must note that, for the particular case $s=3$ Hirschhorn and Sellers proved the identity
$P^{3}_{3}(27n+12)=3P^{3}_{3}(3n+1)$, via generating functions manipulations and some combinatorial arguments [9].
Instead of generating functions, we shall give the formula for $P^{s}_{3}(n)$ (note that $P_{3}(n)=P^{4}_{3}(n)$)
by using representations of posets over the set of natural numbers $\mathbb{N}$, which have been used in [11] to
give criteria for natural numbers which are expressible as sums of three polygonal numbers of positive rank.
\par\bigskip

\textbf{2. Preliminaries}

\par\bigskip

\textbf{2.1. Posets}\par\bigskip

An \textit{ordered set} (or \textit{partially ordered set} or \textit{poset}) is an ordered pair of the form
$(\mathscr{P},\leq)$ of a set $\mathscr{P}$ and a binary relation $\leq$ contained in
$\mathscr{P}\times\mathscr{P}$, called the \textit{order} (or the \textit{partial order}) on $\mathscr{P}$, such
that $\leq$ is reflexive, antisymmetric and transitive [1]. The elements of $\mathscr{P}$ are called the
\textit{points} of the ordered set. \par\bigskip

Let $\mathscr{P}$ be an ordered set and let $x,y\in\mathscr{P}$ we say $x$ is covered by $y$ if $x<y$ and $x\leq
z<y$ implies $z=x$.\par\bigskip

Let $\mathscr{P}$ be a finite ordered set. We can represent $\mathscr{P}$ by a configuration of circles
(representing the elements of $\mathscr{P}$) and interconnecting lines (indicating the covering relation). The
construction goes as follows.

\begin{enumerate}[(1)]
\item To each point $x\in\mathscr{P}$, associate a point $p(x)$ of the Euclidean plane $\mathbb{R}^{2}$, depicted
by a small circle with center at $p(x)$.

\item For each covering pair $x<y$ in $\mathscr{P}$, take a line segment $l(x,y)$ joining the circle  at $p(x)$ to
the circle at $p(y)$.

\item Carry out (1) and (2) in such a way that

\begin{enumerate}[(a)]
\item if $x<y$, then $p(x)$ is lower than $p(y)$,

\item the circle at $p(z)$ does not intersect the line segment $l(x,y)$ if $z\neq x$ and $z\neq y$.

\end{enumerate}

\end{enumerate}
A configuration satisfying (1)-(3) is called a \textit{Hasse diagram} or \textit{diagram} of $\mathscr{P}$. In the
other direction, a diagram may be used to define a finite ordered set; an example is given below, for the ordered
set $\mathscr{P}=\{a,b,c,d,e,f\}$, in which $a< b<c<d<e$, and $f<c$.\par\bigskip

\begin{centering}
\begin{picture}(95,95)

 \setlength{\unitlength}{1pt}
  \setlength{\unitlength}{1pt}

\put(54,30){\line(0,1){20}} \put(54,27){\circle{5}} \put(54,56){\line(0,1){20}} \put(54,53){\circle{5}}
\put(54,79){\circle{5}} \put(52,26){\line(-1,-1){20}} \put(30,4){\circle{5}} \put(26,-12){$f$} \put(63,25){$c$}
\put(63,50){$d$} \put(63,78){$e$} \put(85,2){$b$} \put(111,-22){$a$} \put(78,4){\circle{5}}
\put(80,1){\line(1,-1){20}} \put(56,26){\line(1,-1){20}} \put(102,-21){\circle{5}}

\put(42,-35){Fig. 1}

\end{picture}

\par\bigskip
\end{centering}

\par\bigskip
\par\bigskip
\par\bigskip

 We have only defined diagrams for finite ordered sets. It is not possible to represent the whole of an infinite
ordered set by a diagram, but if its structure is sufficiently regular it can often be suggested
diagrammatically.\par\bigskip

An ordered set $C$ is called a \textit{chain} (or a \textit{totally ordered set} or a \textit{linearly ordered
set}) if and only if  for all $p,q\in C$ we have $p\leq q$ or $q\leq p$ (i.e., $p$ and $q$ are
comparable).\par\bigskip

Let $\mathscr{P}$ be a poset and $S\subset\mathscr{P}$. Then  \textit{$a\in S$ is a maximal element of $S$} if
$a\leq x$ and $x\in S$ imply $a=x$. We denote the set of maximal elements of $S$ by $\mathrm{Max}\hspace{0.1cm}S$.
If $S$ (with the order inherited from $\mathscr{P}$) has a top element, $\top_{S}$ (i.e., $s\leq\top_{S}$ for all
$s\in S$), then $\mathrm{Max}\hspace{0.1cm}S=\{\top_{S}\}$; in this case $\top_{S}$ is called the
\textit{greatest} (or \textit{maximum}) element of $S$, and we write $\top_{S}=\mathrm{max}\hspace{0.1cm}S$.
\par\bigskip

Suppose that $\mathscr{P}_{1}$ and $\mathscr{P}_{2}$ are (disjoint) ordered sets. The \textit{disjoint union}
$\mathscr{P}_{1}+\mathscr{P}_{2}$ of $\mathscr{P}_{1}$ and $\mathscr{P}_{2}$ is the ordered set formed by defining
$x\leq y$ in $\mathscr{P}_{1}+\mathscr{P}_{2}$ if and only if either $x,y\in\mathscr{P}_{1}$ and $x\leq y$ in
$\mathscr{P}_{1}$ or $x,y\in\mathscr{P}_{2}$ and $x\leq y$ in $\mathscr{P}_{2}$. A diagram for
$\mathscr{P}_{1}+\mathscr{P}_{2}$ is formed by placing side by side diagrams for $\mathscr{P}_{1}$ and
$\mathscr{P}_{2}$.\par\bigskip

\textbf{2.2. Partitions and Representations of posets over $\mathbb{N}$}\par\bigskip As usual in this paper
$\mathbb{N}$ denotes the set of natural numbers, while $\mathbb{N}\hspace{0.1cm}\backslash\hspace{0.1cm}\{0\}$ is
the set of positive integers.
\par\bigskip
We denote $t_{k}=p^{3}_{k}=\frac{k(k-1)}{2}$ the $k$-th triangular number $k\in\mathbb{Z}$, and
$s_{k}=p^{4}_{k}=k^{2}$ is the $k$-th square number.\par\bigskip

A \textit{partition} of a positive integer $n$ is a finite nonincreasing sequence of positive integers
$\lambda_{1},\lambda_{2},\dots,\lambda_{r}$ such that $\sum^{r}_{i=1}\lambda_{i}=n$. The $\lambda_{i}$ are called
the \textit{parts} of the partition [6]. A \textit{composition} is a partition in which the order of the summands
is considered.\par\bigskip

Often the partition $\lambda_{1},\lambda_{2},\dots,\lambda_{r}$ will be denoted by $\lambda$ and we sometimes
write $\lambda=(1^{f_{1}}2^{f_{2}}3^{f_{3}}\dots)$ where exactly $f_{i}$ of the $\lambda_{j}$ are equal to $i$.
Note that $\sum_{i\geq1}f_{i}i=n$.\par\bigskip

The \textit{partition function} $p(n)$ is the number of partitions of $n$. Clearly $p(n)=0$ when $n$ is negative
and $p(0)=1$, where the empty sequence forms the only partition of zero. \par\bigskip

Let $(\mathbb{N},\leq)$ be the set of natural numbers endowed with its natural order and $(\mathscr{P},\leq')$ a
poset. A \textit{representation} of $\mathscr{P}$ over $\mathbb{N}$ [11] is a system of the form
\begin{equation}\label{representation}
\begin{split}
\Lambda&=(\Lambda_{0}\hspace{0.1cm};\hspace{0.1cm}(n_{x},\lambda_{x})\mid x\in\mathscr{P}),\\
\end{split}
\end{equation}

where $\Lambda_{0}\subset\mathbb{N}$, $\Lambda_{0}\neq\varnothing$, $n_{x}\in\mathbb{N}$, $\lambda_{x}$ is a
partition with parts in the set $\Lambda_{0}$, and $|\lambda_{x}|$ is the size of the partition $\lambda_{x}$, in
particular if $n_{x}=0$ then we consider $\lambda_{x}=0$. Further

\begin{equation}\label{conditions 1}
\begin{split}
x\leq'y&\Rightarrow n_{x}\leq
n_{y},\quad|\lambda_{x}|\leq|\lambda_{y}|,\quad\text{and}\quad\text{max}\{\lambda_{x}\}\leq\hspace{0.1cm}\text{max}\{\lambda_{y}\}.\\
\end{split}
\end{equation}

\textbf{2.3. The associated graph}\par\bigskip

A \textit{Graph} is a pair $G=(V,E)$ of sets satisfying $E\subseteq V^{2}$, thus the elements of $E$ are
$2$-elements subsets of $V$, such that $V\cap E=\varnothing$. The elements of $V$ are the vertices of the graph
$G$, the elements of $E$ are its edges. A graph with vertex set $V$ is  said to be a graph on $V$. The vertex set
of a graph is referred to as $V(G)$, its edge set as $E(G)$. We write $v\in G$ to a vertex $v\in V(G)$ and $e\in
G$, for an edge $e\in E(G)$, an edge $\{x,y\}$ is usually written as $xy$ or $yx$.\par\bigskip

A vertex $v$ is \textit{incident} with an edge $e$ ; if $v\in e$, then $e$ is an edge at $v$. The two vertices
incident with an edge are its \textit{endvertices} or \textit{ends}, and an edge \textit{joins} its ends. A
\textit{path} is a non-empty graph $P=(V,E)$ of the form\par\bigskip
\begin{centering}

$V=\{x_{0},x_{1},\dots,x_{k}\}$,\quad$E=\{x_{0}x_{1},x_{1}x_{2},\dots,x_{k-1}x_{k}\}$,\par\bigskip

\end{centering}

where the $x_{i}$ are all distinct. The vertices $x_{0}$ and $x_{k}$ are \textit{linked} by $P$ and are called its
\textit{ends}, the vertices $x_{1},\dots,x_{k-1}$ are the inner vertices of $P$. The number of edges of a path is
its \textit{length}, and the path of length $k$ is denoted by $P^{k}$. We often refer to a path by the natural
sequence of its vertices writing $P=x_{0}x_{1}\dots x_{k}=x_{0}||x_{k}$, and calling $P$ a path from $x_{0}$ to
$x_{k}$ [2].\par\bigskip

Given a representation $\Lambda$ for an poset $(\mathscr{P},\leq)$ in [11,12] it was defined its
\textit{associated graph}, $\Gamma_{\Lambda}$ which has as set of vertices the points of $\mathscr{P}$, and
containing all information about partitions of the numbers $n_{x}$. That is $\Gamma_{\Lambda}$ is represented in
such a way that to each vertex of the graph it is attached, either a number $n_{x}$ given by the representation or
one part of a partition of some $n_{y}$ representing some $y\in\mathscr{P}$ such that $x\leq y$.
\par\bigskip

As an example, we consider [11] an infinite sum of infinite chains pairwise incomparable $\mathscr{R}$ in such a
way that $\mathscr{R}=\underset{i=0}{\overset{\infty}{\sum}}C_{i}$, where $C_{j}$ is a chain such that
$C_{j}=v_{0j}<v_{1j}<v_{2j}<\dots$. It is defined a representation over $\mathbb{N}$ for $\mathscr{R}$, by fixing
a number $n\geq3$ and assigning to each $v_{ij}$ the pair
$(n_{ij},\lambda_{ij})=(3+(n-2)i+(n-1)j,(3+(n-2)i+(n-1)j)^{1})$, we note $\mathscr{R}_{n}$ this representation,
and write $v_{ij}\in\mathscr{R}_{n}$ whenever it is assigned the number $n_{ij}=3+(n-2)i+(n-1)j$ to the point
$v_{ij}\in\mathscr{R}$ in this representation. Fig. 2 below suggests the Hasse diagram for this poset with its
associated graph $\Gamma_{p^{n}_{k}}$ which attaches to each vertex $v_{ij}$ the number $3+(n-2)i+(n-1)j$,\quad
$i,j\geq0$.
\par\bigskip
\par\smallskip
 \setlength{\unitlength}{1pt}
  \setlength{\unitlength}{1pt}
  \begin{centering}
\begin{picture}(315,315)

\multiput(40,110)(30,0){10}{\circle{5}} \multiput(40,80)(30,0){10}{\circle{5}}
\multiput(40,50)(30,0){10}{\circle{5}} \multiput(40,20)(30,0){10}{\circle{5}}
\multiput(70,140)(30,0){9}{\circle{5}} \multiput(100,170)(30,0){8}{\circle{5}}
\multiput(130,170)(30,0){6}{\circle{5}}

\multiput(160,200)(30,0){6}{\circle{5}} \multiput(160,230)(30,0){6}{\circle{5}}

\multiput(160,320)(30,0){6}{\circle{5}} \multiput(160,290)(30,0){6}{\circle{5}}
\multiput(160,260)(30,0){6}{\circle{5}}
 \put(40,204){$\Gamma_{p^{n}_{k}}\rightarrow\mathscr{R}$.}
\put(45,20){\vector(1,0){20}} \put(75,20){\vector(1,0){20}}\put(105,20){\vector(1,0){20}}
\put(75,25){\vector(1,1){20}} \put(102,55){\vector(1,2){25}} \put(135,55){\vector(1,1){20}}
\put(165,85){\vector(1,2){25}} \put(132,115){\vector(1,3){26}} \put(105,25){\vector(1,1){20}}
\put(132,55){\vector(1,2){25}} \put(162,115){\vector(1,3){26}} \put(105,50){\vector(1,0){20}}
\put(135,50){\vector(1,0){20}} \put(164,53){\vector(1,1){22}} \put(194,83){\vector(1,1){22}}
\put(194,113){\vector(1,1){22}} \put(194,203){\vector(1,1){22}}  \put(194,143){\vector(1,2){26}}
\put(192,143){\vector(1,3){27}} \put(165,80){\vector(1,0){20}} \put(195,140){\vector(1,0){20}}
\put(195,200){\vector(1,0){20}}

\put(224,143){\vector(1,1){22}} \put(222,143){\vector(1,2){26}} \put(222,204){\vector(1,3){27}}

 \put(225,230){\vector(1,0){20}} \put(225,200){\vector(1,0){20}} \put(224,203){\vector(1,1){22}}  \put(223,234){\vector(1,2){25}}

\put(162,204){\vector(1,4){28}} \put(192,204){\vector(1,4){28}}

\put(165,200){\vector(1,0){20}} \put(255,290){\vector(1,0){20}} \put(284,292){\vector(1,1){22}}
\put(284,262){\vector(1,2){26}}

 \put(282,232){\vector(1,2){26}}
 \put(252,172){\vector(1,2){26}}

 \put(256,204){\vector(1,1){22}} \put(256,234){\vector(1,1){22}}
\put(252,234){\vector(1,2){26}} \put(135,110){\vector(1,0){20}} \put(165,110){\vector(1,0){20}}
\put(165,114){\vector(1,1){22}}
 \put(195,320){\vector(1,0){20}} \put(225,320){\vector(1,0){20}}
 \multiput(40,23)(0,30){3}{\line(0,1){24}}

\multiput(70,23)(0,30){4}{\line(0,1){24}}

\multiput(100,23)(0,30){5}{\line(0,1){24}} \multiput(130,23)(0,30){5}{\line(0,1){24}}
\multiput(160,23)(0,30){10}{\line(0,1){24}} \multiput(190,23)(0,30){10}{\line(0,1){24}}
\multiput(220,23)(0,30){10}{\line(0,1){24}} \multiput(250,23)(0,30){10}{\line(0,1){24}}
\multiput(280,23)(0,30){10}{\line(0,1){24}} \multiput(310,23)(0,30){10}{\line(0,1){24}}
\put(135,55){\vector(1,1){20}}  \put(165,85){\vector(1,2){25}}

\put(224,203){\vector(1,1){22}}

\put(255,174){\tiny{$v_{ij}$}}
\put(26,44){\tiny{$v_{10}$}}\put(56,44){\tiny{$v_{11}$}}\put(56,14){\tiny{$v_{01}$}}
\put(26,14){\tiny{$v_{00}$}}\multiput(40,23)(0,30){3}{\line(0,1){24}}\put(123,186){\tiny{$l.b.p\rightarrow$}}

\put(40,-10){$v_{t_{(i-1)i}}\in l.b.p,\hspace{0.1cm} (\mathrm{left\hspace{0.1cm} boundary\hspace{0.1cm} path})$,}
\put(200,-10){$t_{-1}=0$.}

 \put(165,-35){Fig. 2}

\par\bigskip
\end{picture}
\end{centering}
\par\bigskip
\par\bigskip
\par\bigskip
\par\bigskip
\par\bigskip
\par\bigskip
\par\bigskip
Note that the vertices of the form $v_{(2+7i+3k)(2(i+k)+3)}$, $i,k\geq0$, do not lie on any non-trivial path of
$\Gamma_{p^{n}_{k}}$. \par\bigskip

The representations $\mathscr{R}_{n}$ defined above induce an equivalence relation $\sim$ on $\mathscr{R}$ in such
a way that $v_{ij}\sim v_{kl}$ if and only if $n_{ij}=n_{kl}$. We denote $[v_{ij}]$ the \textit{class} of the
point $v_{ij}\in\mathscr{R}$. Hence the points of\par\bigskip

\begin{centering}

$\mathscr{H}=\underset{k\geq1,m\geq0}{\bigcup}[v_{(2^{k-1}(72^{k}-3)+m(2^{2(k+1)}-1))(32^{k}+2m)}]\cup\underset{i\geq0,s\geq0}{\bigcup}[v_{(2+7i+3s)(2(i+s)+3)}]\subset\mathscr{R}$\par\bigskip

\end{centering}
do not lie on any non-trivial path of $\Gamma_{p^{n}_{k}}$, if $\mathscr{R}$ is represented in such a way that
$(n_{ij},\lambda_{ij})=(2i+j,(2i+j)^{1})$, $i,j\geq0$.\par\bigskip

We say that a path $P\in\Gamma_{p^{n}_{k}}$ is \textit{admissible} if and only if either $P=P_{i0}$, or
$P=P_{i0}P_{ij}$, where $P_{i0}$, $P_{ij}$ are paths such that \par\bigskip

\begin{centering}

$P_{i0}=v_{00}\hspace{0.1cm}||\hspace{0.1cm}v_{t_{i}(i+1)}$,
$i>-1$,\qquad$P_{ij}=v_{t_{i}(i+1)}\hspace{0.1cm}||\hspace{0.1cm}v_{(t_{i}+t_{j})(i+j+2)}$,\par\bigskip

\end{centering}

and the inner vertices of $P_{i0}$ have the form $v_{t_{h}(h+1)}$, $-1<h<i$, while the inner vertices for $P_{ij}$
have the form $v_{(t_{i}+t_{l})(i+l+2)}$, $-1<l<j$. Or $P=P_{i0}P_{ij}P_{ijk}$, where $P_{ijk}$ is a path such
that\par\bigskip
\begin{centering}

$P_{ijk}=v_{(t_{i}+t_{j})(i+j+2)}\hspace{0.1cm}||\hspace{0.1cm}v_{(t_{i}+t_{j}+t_{k})(i+j+k+3)}$, $-1<k\leq
j$.\par\bigskip

\end{centering}

If for $s\geq3$ fixed, we consider the representation $\mathscr{R}_{s}$ then the numbers $n_{t_{i(i+1)}}$
representing the vertices $v_{t_{i(i+1)}}\in\mathscr{R}_{s}$ of the \textit{left boundary path}, $l.b.p$ are
expressible in the form $n_{t_{i}(i+1)}=p^{s}_{1}+p^{s}_{1}+p^{s}_{i+2}$, $s$ is an fixed index, $i\geq-1$.
Furthermore if $i_{0}\geq0$ is a fixed number, and $v_{t_{i_{0}(i_{0}+1)}}\in l.b.p$\hspace{0.1cm}
then\hspace{0.1cm} $n_{(t_{i_{0}}+t_{l})(i_{0}+l+2)}$ representing the vertex $v_{(t_{i_{0}}+t_{l})(i_{0}+l+2)}\in
P_{i_{0}}P_{ij}$, can be expressible in the form
$n_{(t_{i_{0}}+t_{l})(i_{0}+l+2)}=p^{s}_{i_{0}+2}+p^{s}_{l+2}+p^{s}_{1}$, and a number
$n_{(t_{i_{0}}+t_{j}+t_{k})(i_{0}+j+k+3)}$ has the form
$n_{(t_{i_{0}}+t_{j}+t_{k})(i_{0}+j+k+3)}=p^{s}_{i_{0}+2}+p^{s}_{j+2}+p^{s}_{k+2}$. In [11] has been used this
facts to state the following theorem ;

\begin{teor}\label{component}
A number $m\in\mathbb{N}\hspace{0.1cm}\backslash\hspace{0.1cm}\{0\}$ is the sum of three $n$-gonal numbers of
positive rank if and only if $m$ represents a vertex $v_{ij}\in\mathscr{R}_{n}$ in a non-trivial component of
$\Gamma_{p^{n}_{k}}$.

\end{teor}
For example, if $\mathscr{R}'_{3}$ is a representation such that $(n_{ij},\lambda_{ij})=(i+j,(i+j)^{1})$,
$i,j\geq0$, and $\mathscr{G}_{i}$ is the family of subsets of $\mathscr{R}$ such that for $i\geq0$,
$\mathscr{G}_{i}=\{[v_{(t_{i}+j)(i+2)}]\mid0\leq j\leq t_{i+1}\}$ then
$V(\Gamma_{p^{n}_{k}})=\underset{i}{\bigcup}(\mathscr{G}_{i}\cup [v_{00}]\cup[v_{01}])$, because every natural
number is expressible as a sum of three or fewer triangular numbers.\par\bigskip

\textbf{3. The main result}\par\bigskip

If $v_{i_{0}j}\in\Gamma_{p^{n}_{k}}$ belongs to an non-trivial component, and $v_{i_{0}j}\in\mathscr{R}_{s}$ then
there exists an admissible path $P=v_{00}\hspace{0.1cm}||\hspace{0.1cm}v_{i_{0}j}$, which has associated a family
of compositions (see 2.2) in such a way that if $l>0$ then
\begin{equation}\label{partitions}
\begin{split}
n_{t_{i_{0}}(i_{0}+1)}&=p^{s}_{i_{0}+2}+p^{s}_{1}+p^{s}_{1},\hspace{0.2cm}i_{0}\geq-1\hspace{0.2cm}\mathrm{if}\hspace{0.1cm}P=P_{i_{0}0},\\
n_{(t_{i_{0}}+t_{l})(i_{0}+l+2)}&=p^{s}_{l+2}+p^{s}_{i_{0}+2}+p^{s}_{1}=p^{s}_{i_{0}+2}+p^{s}_{l+2}+p^{s}_{1},\hspace{0.1cm}\mathrm{if}\hspace{0.1cm}P=P_{i_{0}0}P_{i_{0}j},\raisetag{15mm}\\
n_{(t_{i_{0}}+t_{j}+t_{k})(i_{0}+j+k+3)}&=p^{s}_{i_{0}+2}+p^{s}_{j+2}+p^{s}_{k+2}=p^{s}_{k+2}+p^{s}_{j+2}+p^{s}_{i_{0}+2},\hspace{0.2cm}\mathrm{and}\\
n_{(t_{i_{0}}+t_{j}+t_{k})(i_{0}+j+k+3)}&=p^{s}_{j+2}+p^{s}_{k+2}+p^{s}_{i_{0}+2},\hspace{0.1cm}-1\leq i_{0}\leq k\leq j,\hspace{0.1cm}\mathrm{if}\hspace{0.1cm}P=P_{i_{0}0}P_{i_{0}j}P_{i_{0}jk}.\\
\end{split}
\end{equation}

Thus we say that two admissible paths $P$, $Q$ are \textit{equivalent} if and only if they have associated the
same partitions. If $P=v_{00}\hspace{0.1cm}||\hspace{0.1cm}v_{ij}$ is an admissible path then we note $[P]$ the
class of $P$. Therefore $[P]$ can be determined by fixing the end $v_{ij}$, we note $\overline{v_{ij}(P)}$
whenever a vertex $v_{ij}$ has been fixed due to this condition (i.e, $\overline{v_{ij}(P)}= [P]$,
$\overline{w_{ij}(Q)}=\varnothing$ for $Q\in[P]$, $w_{ij}\in[v_{ij}]$, and $w_{ij}\neq v_{ij}$). For example
$\overline{v_{t_{i}(i+1)}(P_{t_{i}0})}=[P_{t_{i}0}]$ only contains the admissible path $P_{t_{i}0}$.\par\bigskip

If $v_{ij}(P)$ is the set of admissible paths of the form $v_{00}||u_{rs}$ where $u_{rs}\in[v_{ij}]$ then we
note\par\bigskip

$\mathscr{A}(v_{ij})=\{w_{ij}\in[v_{ij}]\mid\overline{w_{ij}(P)}=[P], P\in v_{ij}(P)\}$, therefore
$v_{ij}(P)=\underset{w_{ij}\in\mathscr{A}(v_{ij})}{\bigcup}\overline{w_{ij}(P)}$,\par\bigskip

If $\mathscr{A}(v_{ij})\neq\varnothing$ and $w_{ij}\in\mathscr{A}(v_{ij})$ then $\delta_{g}(w_{ij})$ denotes the
number of classes of admissible paths at $w_{ij}$. If $\mathscr{A}(v_{ij})=\varnothing$ then
$\delta_{g}(v_{ij})=0$.\par\bigskip The next theorem is a consequence of theorem \ref{component} and the equations
(\ref{partitions}).

\begin{teor}\label{three pentagonal}
Let $P^{s}_{3}(n)$ denote the number of partitions of $n\in\mathbb{N}$ into three $s$-gonal numbers, $s\geq3$.
then
\[P^{s}_{3}(n)=
\begin{cases}
\underset{w_{ij}\in\mathscr{A}(v_{ij})}{\sum}\delta_{g}(w_{ij}), &\mathrm{if}\hspace{0.2cm}\mathscr{A}(v_{ij})\neq\varnothing,\\
\hspace{1.2cm}                    0, &\mathrm{otherwise}.
\end{cases}\]

where $v_{ij}\in\mathscr{R}_{s}$ and $n_{ij}=n$.\hspace{0.5cm}\text{\qed}

\end{teor}

For example in $\mathscr{R}_{4}$, $[v_{(10,5)}]\cap\mathscr{A}(v_{(10,5)})=\{v_{(10,5)},v_{(7,7)}\}$,
$\delta_{g}(v_{(10,5)})=1$, $\delta_{g}(v_{(7,7)})=1$,\par\bigskip If
$P=v_{(0,0)}||v_{(0,1)}||v_{(6,5)}||v_{(7,7)}$,\quad $Q=v_{(0,0)}||v_{(10,5)}$ then\par\bigskip

$\overline{v_{(7,7)}(P)}=\{v_{(0,0)}||v_{(0,1)}||v_{(6,5)}||v_{(7,7)},v_{(0,0)}||v_{(6,4)}||v_{(7,6)}||v_{(7,7)},v_{(0,0)}||v_{(2,1)}||v_{(7,6)}||v_{(7,7)}\}$.\par\bigskip
$\overline{v_{(10,5)}(Q)}=\{v_{(0,0)}||v_{(10,5)}\}$. Therefore $P^{4}_{3}(38)=2$.

\end{document}